\documentclass[12pt,a4paper]{article}
\usepackage{amsmath,amssymb,amsfonts}
\usepackage[latin1]{inputenc}

\def\Bbb{\mathbb}

\title{\bf On the decimal and octal digits of $1/p$}

\author{Kurt Girstmair}

\date{}

\makeatletter
\let\@@maketitle=\maketitle
\def\maketitle{\def\thispagestyle##1{\relax}\@@maketitle}
\makeatother
\newtheorem{theorem}{Theorem}

\newtheorem{lemma}{Lemma}

\def\BE{\begin{equation}}
\def\EE{\end{equation}}
\def\BD{\begin{displaymath}}
\def\ED{\end{displaymath}}
\def\BA{\begin{array}}
\def\EA{\end{array}}
\def\BEA{\begin{eqnarray*}}
\def\EEA{\end{eqnarray*}}
\def\BI{\bibitem}

\def\Z{\Bbb Z}
\def\Q{\Bbb Q}

\def\phi{\varphi}

\def\MB{\mbox}
\def\LD{\ldots}
\def\OV{\overline}
\def\SP#1{\langle #1 \rangle}

\def\WT{\widetilde}

\def\sminus{\smallsetminus}

\def\NDIV{\, \nmid \,}

\def\MN{\medskip\noindent}
\def\STOP{\hfill$\Box$}

\def\LS#1#2{ \left( \frac{#1}{#2} \right) }
\def\LF#1#2{\left\lfloor\frac{#1}{#2}\right\rfloor}

\begin{document}
\maketitle

\begin{abstract}

\noindent
Let $p$ be a prime $\equiv 3$ mod 4, $p>3$, and suppose that 10 has the order $(p-1)/2$ mod p. Then $1/p$ has a decimal period of length
$(p-1)/2$. We express the frequency of each digit $0,\LD,9$ in this period in terms of the class numbers of two imaginary quadratic number fields.
We also exhibit certain analogues of this result, so for the case that 10 is a primitive root mod $p$ and for the octal digits of $1/p$.

\end{abstract}

%%%%%%%%%%%%%%%%%%%%%%%%%%%%%%%%%%%%%%%%%%%%%
\section*{1. Introduction}
%%%%%%%%%%%%%%%%%%%%%%%%%%%%%%%%%%%%%%%%%%%%%

A connection between the digits of $1/p$ and class numbers was first established in \cite{Gi1}. A number of papers on this topic appeared in the sequel; see \cite{Gi2}, \cite{Hi}, \cite{MuTh},
\cite{ChKr},  \cite{Mi}, \cite{PuSa}, \cite{To}, and \cite{Sh}. In the present article, the focus lies on the decimal digits of $1/p$ or, more generally, of $m/p$, where $m\in\{1,\LD,p-1\}$ is a quadratic residue mod $p$.

Let $p$ be a prime, $p\equiv 3$ mod 4, $p>3$. Suppose, moreover, that $10$ has the order $(p-1)/2$ mod $p$. Then $m/p$ has the decimal expansion
\BD
\label{1.1}
 m/p=\sum_{j=1}^{\infty}a_j10^{-j},
\ED
where the numbers $a_j\in\{0,\LD,9\}$ are the digits of $m/p$. The sequence $(a_1,\LD, a_{(p-1)/2})$ is the period of this expansion.

The connection with class numbers in this case, as known so far, does not give much information about the distribution of the digits in the period. Typically, most of the aforementioned authors determine
certain mean values of the digits, for instance,
\BD
  \sum_{j=1}^{(p-1)/2}a_j=9(p-1)/4-9h_1/2,
\ED
where $h_1$ is the class number of the imaginary quadratic number field $\Q(\sqrt{-p})$; see \cite[Satz 11]{Gi2}.

In this paper, however, we exhibit the frequency of each digit in the period.
For $k\in\{0,\LD,9\}$ we define
\BE
\label{1.1.1}
  n_k=|\{j; 1\le j\le (p-1)/2, a_j=k\}|.
\EE
Let $h_1$ be as above and $h_2$ the class number of $\Q(\sqrt{-5p})$.
The frequencies $n_k$, $k=0,\LD,9$, can be expressed in terms of $h_1$ and $h_2$.

\begin{theorem} % Theorem 1 %%%%%%%%%%%%%%%%%%%%%%%%%%%%%%%%%%%%%%%%%%%%%%
\label{t1}

In the above setting,
\BD
  n_k=\frac 12\left(\LF{(k+1)p}{10}-\LF{kp}{10}+\delta_k\right) \MB{ and }
  n_{9-k}=n_k-\delta_k, \enspace k=0,\LD, 4,
\ED
where the numbers $\delta_k$ take the following values. If $p\equiv 3$ mod 8, then
\BD
  \delta_0 = 0,
  \delta_1 = \delta_4=3h_1/2-h_2/4,
  \delta_2 = -\delta_1,
  \delta_3 = 3h_1/2+h_2/4.
\ED
If $p\equiv 7$ mod 8, then
\BD
  \delta_0 = h_1-h_2/2,
  \delta_1 = h_2/4,
  \delta_2 = 3h_2/4,
  \delta_3 =\delta_4=-\delta_1.
\ED

\end{theorem} %%%%%%%%%%%%%%%%%%%%%%%%%%%%%%%%%%%%%%%%%%%%%%%%%%%%%%%%%%%%

\MN
The proof of this theorem is based on a result of B. C. Berndt; see Section 2. The present author's contribution is of a more elementary nature, but it gives a result that is worth knowing.
It seems that no result about the digits that involves two class numbers can be found in the literature.

The first six primes $p$ falling under Theorem \ref{t1} are 31, 43, 67, 71, 83, 107.

\MN
{\em Example.} Let $p=67$. Then $h_1=1$, $h_2=18$. We have
\BE
\label{1.2}
 1/67=0.\OV{014925373134328358208955223880597},
\EE
where the bar marks the period. From Theorem 1 we obtain $\delta_0=0$, $\delta_1=-3$, $\delta_2=3$, $\delta_3=6$, $\delta_4=-3$, and
$n_0=3$, $n_1=2$, $n_2=5$, $n_3=6$, $n_4=2$, $n_5=5$, $n_6=0$, $n_7=2$, $n_8=5$, $n_9=3$, in accordance with (\ref{1.2}).

\MN
{\em Remarks.} 1. We have $p\equiv l$ mod 10 for $l\in\{1,3,7,9\}$. Then
\BE
\label{1.1.2}
   \LF{(k+1)p}{10}-\LF{kp}{10}=\LF{p}{10},
\EE
except for
\BEA
l&=&3 \MB{ and } k=3,\\
l&=&7 \MB{ and } k=1,2,4,\\
l&=&9 \MB{ and }k=1,2,3,4.
\EEA
In these cases the left-hand side of (\ref{1.1.2}) equals $\LF{p}{10}+1$.

2. Some consequences of Theorem \ref{t1} are immediate. Indeed, if $p\equiv 3$ mod 8, then $n_3>n_6$. If $p\equiv 7$ mod 8, then $n_1>n_8$, $n_2>n_7$, $n_3<n_6$, $n_4<n_5$.

3. Under the Generalized Riemann Hypothesis, the natural density of the primes $p\equiv 3$ mod 4 such that 10 has the order $(p-1)/2$ mod p is $A/2=0.186977\LD,$ where $A$ is Artin's constant
\BE
\label{1.30}
   A=\prod_{q}\left(1-\frac 1{q(q-1)}\right)=0.3739558\LD,
\EE
$q$ running through all primes. The author was informed about this fact by P. Moree, who also pointed out that
\cite[p. 663]{MoZu} contains an analogous computation for the number 3 instead of 10.
Hence one expects that about 37\% of all primes $p\equiv 3$ mod 4 below a large bound have this property. If we include the primes considered in Theorem \ref{t2} below,
we have even more than 74\% of all primes $p\equiv 3$ mod 4 below such a bound.

4. There are fairly effective methods for computing the class numbers of imaginary quadratic number fields, see \cite{Lo}. Therefore, the frequencies $n_0,\LD,n_9$ of the digits
can be computed by Theorem \ref{t1} for primes $p$ of an order of magnitude like $10^{15}$, where a naive computation is hopelessly slow.

5. Of course, one may ask what happens in the case $p\equiv 1$ mod 4. Formula (6) in \cite{Sz} yields a result for this case. Suppose that $10$ has the order $(p-1)/2$ mod p and $n_k$ is defined by
(\ref{1.1.1}). We obtain, for instance,
$n_0=(\lfloor p/10\rfloor +\delta_0)/2$, with
\BE
\label{1.20}
  \delta_0=\frac 12\MB{ Re\,}\left(\OV{\psi}(-p)(1-\chi\psi(2))B_{\chi\psi}\right),
\EE
where $\psi$ is the (odd) Dirichlet character mod 5 defined by $\psi(2)=i$, $\OV{\psi}$ the complex-conjugate character, $\chi$ the Legendre symbol $\LS{}{p}$,
$\chi\psi$ the character mod $5p$ defined by $\chi\psi(j)=\chi(j)\psi(j)$
and
\BD
  B_{\chi\psi}=\frac 1{5p}\sum_{j=1}^{5p}j \chi\psi(j)
\ED
the corresponding Bernoulli number of order 1. However, we see no way to an interpretation of (\ref{1.20}) in terms of class numbers.

\MN
In Section 2 we prove Theorem \ref{t1}
and show what is possible if the order of 10 mod $p$ is an arbitrary odd number.
We also give an analogue of Theorem \ref{t1} for the case that 10 is a primitive root mod $p$. In Section 3 we note the analogue of this theorem for the octal expansion of $m/p$.
In this case, however, we can also treat primes $p\equiv 1$ mod 4.

%%%%%%%%%%%%%%%%%%%%%%%%%%%%%%%%%%%%%%%%%%%%%
\section*{2. Digits and quadratic residues}
%%%%%%%%%%%%%%%%%%%%%%%%%%%%%%%%%%%%%%%%%%%%%

Let $p$ be a prime, $b\ge 2$ an integer with $p\NDIV b$, and $m\in\{1,\LD,p-1\}$. Then $m/p$ has the expansion
\BE
\label{2.2}
  m/p=\sum_{j=1}^{\infty}a_jb^{-j},
\EE
where the numbers $a_j\in\{0,\LD,b-1\}$ are the digits of $m/p$ with respect to the basis $b$. For an integer $k$ let $(k)_p$ denote the number $j\in\{0,\LD,p-1\}$ that satisfies $k\equiv j$ mod p.
We define
\BD
  \theta_b(k)=\frac{b(k)_p-(bk)_p}p
\ED
for $k\in \Z$. Then
\BE
\label{2.3}
  a_j=\theta_b(mb^{j-1}),\enspace j\ge 1,
\EE
see \cite{Gi2}. This shows, in particular, that $(a_1,\LD,a_q)$, where $q$ is the order of $b$ mod p, is a period of the expansion (\ref{2.2}).

The basic idea of Theorem \ref{t1} consists in establishing  a connection between the said digits and integers in certain subintervals of $[0,p]$.
This connection is contained in the following lemma.

\begin{lemma} % Lemma 1 %%%%%%%%%%%%%%%%%%%%%%%%%%%%%%%%%%%%%%%%%%%%%%%%
\label{l1}
Let $p$ be a prime, $b\ge 2$, and $l\in\{1,\LD,p-1\}$. Let $k$ be an integer, $0\le k\le b-1$. Then
\BD
   \theta_b(l)=k
\ED
if, and only if,
\BD
 \frac{kp}b<l<\frac{(k+1)p}b.
\ED
\end{lemma} %%%%%%%%%%%%%%%%%%%%%%%%%%%%%%%%%%%%%%%%%%%%%%%%%%%%%%%%%%%%

\MN{\em Proof.}
If $\theta_b(l)\le k$, then $l \le (kp+(bl)_p)/b\le (kp+p-1)/b$ and $l<(k+1)p/b$.
If $\theta_b(l)\ge k$, then $l\ge (kp+(bl)_p)/b$ and $l>kp/b$.

Conversely, if $l>kp/b$, then $\theta_b(l)>k-(bl)_p/p>k-1$. Since $\theta_b(l)$ is an integer, we have $\theta_b(l)\ge k$.
If $l<(k+1)p/b$, then $\theta_b(l)<k+1-(bl)_p/p<k+1$ and $\theta_b(l)\le k$.
\STOP

\MN
We collect further ingredients of the proof of Theorem \ref{t1}.
By $Q$  we denote the set of quadratic residues in $\{1,\LD,p-1\}$ and by $N$ the set $\{1,\LD,p-1\}\sminus Q$.
Let $b\ge 2$ be as above. Then
\BE
\label{2.4}
   \left|\Z\cap \left(\frac{kp}b,\frac{(k+1)p}b\right)\right|=\LF{(k+1)p}{b}-\LF{kp}{b}.
\EE
Moreover, let $p\equiv 3$ mod 4. For $k\in\{0,\LD,b-1\}$ put $k'=b-1-k$. Then
\BE
\label{2.6}
   \left|Q\cap \left(\frac{k'p}b,\frac{(k'+1)p}b\right)\right|= \left|N\cap \left(\frac{kp}b,\frac{(k+1)p}b\right)\right|.
\EE
Indeed, an integer $l$ lies in $Q$ if, and only if, $p-l$ lies in $N$.

The following lemma is a special case of Theorem 8.1 in \cite{Be}, which, however, must be expressed in terms of class numbers; for this purpose
see \cite[Cor. 4.6, Th. 4.9, Th. 4.17]{Wa} and \cite[p. 68]{Ha}. This lemma involves the Legendre symbol $\LS{\:}{p}$, which we denote by $\chi$ for typographical reasons.

\begin{lemma} % Lemma 2 %%%%%%%%%%%%%%%%%%%%%%%%%%%%%%%%%%%%%%%%%%%%%%%%
\label{l2}

Let $p\equiv 3$ mod 4, $p>3$.
For $k=0,\LD,4$,
\BD
  \left|Q\cap\left(\frac{kp}{10},\frac{(k+1)p}{10}\right)\right|-\left|N\cap\left(\frac{kp}{10},\frac{(k+1)p}{10}\right)\right|=\delta_k
\ED
with
\BEA
  \delta_0&=& (3+\chi(2)+\chi(5)-\chi(10))h_{-p}/4-(1+\chi(2))h_{-5p}/4,\\
  \delta_1&=&(2-\chi(2))(1-\chi(5))h_{-p}/4+\chi(2)h_{-5p}/4,\\
  \delta_2&=&(2-\chi(2))(\chi(5)-1)h_{-p}/4+(2+\chi(2))h_{-5p}/4,\\
  \delta_3&=&(2-\chi(2))(1-\chi(5))h_{-p}/4-\chi(2)h_{-5p}/4,\\
  \delta_4&=&(3-4\chi(2)+\chi(5))h_{-p}/4-h_{-5p}/4.
\EEA

\end{lemma} %%%%%%%%%%%%%%%%%%%%%%%%%%%%%%%%%%%%%%%%%%%%%%%%%%%%%%%%%%%%

\MN
{\em Proof of Theorem} \ref{t1}. Let $p\equiv 3$ mod 4, $m\in Q$, and $(p-1)/2$ be the order of 10 mod $p$.  We consider the expansion (\ref{2.2}) in the case $b=10$.
Observe that the numbers $(m10^{j-1})_p$, $j=1,\LD,(p-1)/2$, run through $Q$.
Let $k\in\{0,\LD,4\}$.
From (\ref{2.3}) and Lemma \ref{l1} we see that
\BE
\label{2.14}
  n_k=\left|Q\cap \left(\frac{kp}{10},\frac{(k+1)p}{10}\right)\right|.
\EE
Now (\ref{2.4}) gives
\BD
  \left|Q\cap \left(\frac{kp}{10},\frac{(k+1)p}{10}\right)\right|+\left|N\cap \left(\frac{kp}{10},\frac{(k+1){p}}{10}\right)\right|=\LF{(k+1)p}{10}-\LF{kp}{10}.
\ED
Combined with Lemma \ref{l2}, this yields the first assertion of Theorem \ref{t1}, however, with $\delta_k$ as in Lemma \ref{l2}.
But since $\chi(10)=1$, only two cases must be distinguished. If $p\equiv 3$ mod 8, then $\chi(2)=-1$. Accordingly, $\chi(5)=-1$, and the values of the numbers $\delta_k$ are as in the
respective case of Theorem \ref{t1}. If $p\equiv 7$ mod 8, then $\chi(2)=\chi(5)=1$, and we obtain the values of $\delta_k$ as given in Theorem \ref{t1} for this case.
The second assertion follows from (\ref{2.6}).
\STOP

\MN
{\em Remarks.} 1. In the case $m\in N$, Theorem \ref{t1} remains valid, provided that the respective signs of the numbers $\delta_k$ are interchanged.
The proof is a simple variation of the above proof.

2. One can derive Lemma \ref{l2} from Formula (6) in \cite{Sz}, but this is somewhat laborious. In fact, this formula settles not only the case $p\equiv 1$ mod 4,
as we mentioned in the Introduction, but also our case $p\equiv 3$ mod 4.

\MN
Assume, in the above setting, that $10\in Q$ has the order $q$ mod $p$. Then $q$ is a divisor of the odd number $(p-1)/2$.
Let $H$ be the group of squares in $G=(\Z/p\Z)^{\times}$ and $b_1,\LD b_{(p-1)/(2q)}$ a system of representatives of the group $H/\SP{\OV{10}}$ in the set $\{1,\LD,p-1\}$.
Then
\BD
  b_l/p=\sum_{j=1}^{\infty}a_j^{(l)}10^{-j}, \enspace l=1,\LD,(p-1)/(2q)
\ED
with $a_j^{(l)}=\theta_b(b_l10^{j-1})$, $j\ge 1$. The fraction $b_l/p$ has the period $(a_1^{(l)},\LD,a_q^{(l)})$.
For $k=0,\LD,9$ we consider
\BE
\label{2.12}
  \WT n_k=|\{(l,j); a_j^{(l)}=k\}|,
\EE
where $l$ and $j$ run through $1,\LD,(p-1)/2q$ and $1,\LD,q$, respectively. It is easy to see that the assertions of Theorem \ref{t1} remain valid if
$n_k$ is replaced by $\WT n_k$, $k=0,\LD,9$.

\MN
{\em Example.} Let $p=79$. Then $10$ is a quadratic residue mod $p$ of the order $q=13$. So $(p-1)/(2q)=3$.
We use the primitive root 3 mod $p$ and obtain $b_1=1$, $b_2=(3^2)_p=9$, $b_3=(3^4)_p=2$.
Now
\BD
 1/79=0.\OV{0126582278481},\: 9/79=0.\OV{1139240506329},\: 2/79=0.\OV{0253164556962}.
\ED
Here $p\equiv 7$ mod 8, $h_1=5$, and $h_2=8$. Theorem \ref{t1} yields $\delta_0=1,\delta_1=2$, $\delta_2=6$, $\delta_3=-2=\delta_4$.
This gives $\WT n_0=4$, $\WT n_1=5$, $\WT n_2=7$, $\WT n_3=3=\WT n_4$, $\WT n_5=5=\WT n_6$, $\WT n_7=1$, $\WT n_8=3=\WT n_9$,
which agrees with the displayed periods.

\MN
Let $p$ be a prime $\equiv 3$ mod 4, $p>3$.
Let $10\ge 2$ be  a primitive root mod $p$.
Moreover, let $m\in Q$.
Then $m/p$ has the expansion (\ref{2.2}) with $a_j=\theta_b(mb^{j-1})$; recall (\ref{2.3}).
Accordingly, $(a_1,\LD,a_{p-1})$ is a period of (\ref{2.2}).
We define
\BD
  n_k=|\{j; j \MB{ odd, }a_j=k\}| \: \MB{ and }\: n_k^*=|\{j; j \MB{ even, }a_j=k\}|,
\ED
where $j$ runs through $1,\LD,p-1$.
Now we have, instead of (\ref{2.14}),
\BD
   n_k=\left|Q\cap \left(\frac{kp}{10},\frac{(k+1)p}{10}\right)\right| \:\MB{ and }\:
   n_k^*=\left|N\cap \left(\frac{kp}{10},\frac{(k+1)p}{10}\right)\right|.
\ED
Proceeding as in the proof of Theorem \ref{t1}, we obtain an analogous result. Here, however, $\chi(10)=-1$, which changes the values of the numbers $\delta_k$. Indeed, for
$p\equiv 3$ mod 8 we obtain
\BE
\label{2.16}
  \delta_0 = h_1,
  \delta_1 = -h_2/4,
  \delta_2 = \delta_3=-\delta_1,
  \delta_4 =2h_1-h_2/4.
\EE
If $p\equiv 7$ mod 8, then
\BE
\label{2.18}
  \delta_0 = h_1-h_2/2,
  \delta_1 = h_1/2+h_2/4,
  \delta_2 = -h_1/2+3h_2/4,
  \delta_3 = h_1/2-h_2/4,
  \delta_4 = -\delta_1.
\EE

\begin{theorem} % Theorem 2 %%%%%%%%%%%%%%%%%%%%%%%%%%%%%%%%%%%%%%%%%%%%%%
\label{t2}
In the above setting, we have,
for $k=0,\LD,4$,
\BD
\label{3.2}
  n_k=\frac 12\left(\LF{(k+1)p}{10}-\LF{kp}{10}+\delta_k\right)\: \MB{ and }\:
  n_{9-k}=n_k-\delta_k,
\ED
where $\delta_k$ is as in {\rm (\ref{2.16})} if $p\equiv 3$ mod 8, and as in {\rm (\ref{2.18})} if $p\equiv 7$ mod 8. Moreover, $n_k^*=n_{9-k}$ for $k=0,\LD,9$.

\end{theorem} %%%%%%%%%%%%%%%%%%%%%%%%%%%%%%%%%%%%%%%%%%%%%%%%%%%%%%%%%%%%

\MN
{\em Example.} Let $p=47$. Then 10 has the order 46 mod $p$ and $p\equiv 7$ mod 8. Here $h_1=5$ and $h_2=2$.
From (\ref{2.18}) we obtain $\delta_0=4$, $\delta_1=3$, $\delta_2=-1$, $\delta_3=2$, $\delta_4=-3$. The values of $n_k$ are 4, 4, 2, 3, 1, 4, 1, 3, 1, 0 in the order corresponding to $k=0,\LD,9$.
The values of $n_k^*$ are the same, but in the converse order.

\MN
{\em Remark.} Under the Generalized Riemann Hypothesis, the natural density of primes $p\equiv 3$ mod 4 such that 10 is a primitive root mod p is $A/2$, where $A$ is Artin's constant (recall
(\ref{1.30})); see Theorem 2 in \cite{Mo}. Hence one expects that more than 74\% of all primes $\equiv 3$ mod 4 below a large bound fall under Theorem \ref{t1} or Theorem \ref{t2}.

%%%%%%%%%%%%%%%%%%%%%%%%%%%%%%%%%%%%%%%%%%%%%
\section*{3. Octal digits}
%%%%%%%%%%%%%%%%%%%%%%%%%%%%%%%%%%%%%%%%%%%%%

In this section we give an analogue of Theorem \ref{t1} for the octal expansion of $1/p$, both for $p\equiv 3$ mod 4 and $p\equiv 1$ mod 4. First let $p\equiv $ 3 mod 4, $p>3$,
and suppose that $b=8$ has the order $(p-1)/2$ mod p.
Accordingly, $\chi(8)=\chi(2)=1$ and $p\equiv 7$ mod 8.
Then $1/p$ has the expansion (\ref{2.2}) with $b=8$ and $a_j\in\{0,\LD,7\}$. The period is $(a_1,\LD,a_{(p-1)/2})$. For $k\in\{0,\LD,7\}$ let $n_k$ denote the frequency of the digit $k$ in this period, i.e.,
$n_k$ is defined as in (\ref{1.1.1}). By $h_1$ we denote the class number of $\Q(\sqrt{-p})$ and by $h_2$ the class number of $\Q(\sqrt{-2p})$.

\begin{theorem} % Theorem 3 %%%%%%%%%%%%%%%%%%%%%%%%%%%%%%%%%%%%%%%%%%%%%%
\label{t3}

In the above setting,
\BD
  n_k=\begin{cases}
      \frac 12\left(\LF{p}{8}+\delta_0\right)   &\MB{ if } k=0;\\
      \frac 12\left(\LF{p}{8}+1+\delta_k\right) &\MB{ if } k\in\{1,2,3\},
      \end{cases}
\ED
where
$
  \delta_0= h_1-h_2/4,
  \delta_1=\delta_2=h_2/4,
  \delta_3=-\delta_1,
$
Moreover, $n_{7-k}=n_k-\delta_k$, $k=0,\LD, 3$.

\end{theorem} %%%%%%%%%%%%%%%%%%%%%%%%%%%%%%%%%%%%%%%%%%%%%%%%%%%%%%%%%%%%

\MN
The proof of this theorem  follows the above pattern. The crucial ingredient is Theorem 7.1 of \cite{Be}, which has to be interpreted by means of \cite[Cor. 4.6, Th. 4.9, Th. 4.17]{Wa}
and \cite[p. 68]{Ha}.

Now suppose $p\equiv 1$ mod 4 and let $(p-1)/2$ be the order of $8$ mod p. The number 2 is a quadratic residue mod $p$,
whence $p\equiv 1$ mod 8 follows.
Let $h_1$ be the class number of $\Q(\sqrt{-p})$ and $h_2$ the class number of $\Q(\sqrt{-2p})$.

\begin{theorem} % Theorem 4 %%%%%%%%%%%%%%%%%%%%%%%%%%%%%%%%%%%%%%%%%%%%%%
\label{t4}

In the above setting,
\BD
  n_k=
      \frac 12\left(\LF{p}{8}+\delta_k\right),
\ED
$k=0,\LD,3$, where
$
  \delta_0= (h_1+h_2)/4,
  \delta_1=\delta_3=(h_1-h_2)/4,
  \delta_2=(-3h_1+h_2)/4.
$
Moreover, $n_{7-k}=n_k$, $k=0,\LD,3$.

\end{theorem} %%%%%%%%%%%%%%%%%%%%%%%%%%%%%%%%%%%%%%%%%%%%%%%%%%%%%%%%%%%%

\MN
For the proof we use Theorem 7.1 of \cite{Be} again and the symmetry
\BD
 \left(\frac{kp}8,\frac{(k+1)p}8\right)\cap Q=\left(\frac{k'p}8,\frac{(k'+1)p}8\right)\cap Q,
\ED
with $k'=7-k$, $k=0,\LD,3$. This symmetry is due to $-1\in Q$.

\bigskip
\centerline {\bf Acknowledgment}

\MN
The author expresses his gratitude to Pieter Moree for his information about the density of primes $p\equiv 3$ mod 4 such that $10$ has the order $(p-1)/2$ mod $p$.

\bigskip
\centerline{\bf Declaration}

\noindent
This preprint has not undergone
peer review or any post-submission improvements or corrections. The Version of Record of this article is
published in {\em Archiv der Mathematik}, and is available online at https://doi.org/10.1007/s00013-026-02252-z.

%%%%%%%%%%%%%%%%%%%%%%%%%%%%%%%%%%%%%%%%%%%%%%%%%%%%%
%%%%%%%%%%%%%%%%%%%%%%%%%%%%%%%%%%%%%%%%%%%%%%%%%%%%%%%%%%%%%%%%%%%%%%%%%%

\MN
Kurt Girstmair\\
Institut f\"ur Mathematik \\
Universit\"at Innsbruck   \\
Technikerstr. 13/7        \\
A-6020 Innsbruck, Austria \\
Kurt.Girstmair@uibk.ac.at

\end{document}